\documentclass[11pt,openbib]{article}
\usepackage{amsmath}
\usepackage{amssymb}
\usepackage{amsfonts}
\usepackage[T1]{fontenc}
\usepackage{tikz-cd}

\newcommand{\R}{\mathbb{R}}

\newcommand{\Z}{\mathbb{Z}}
\newcommand{\Gl}{\mathrm{Gl}}
\newcommand{\gl}{\mathrm{gl}}

\newcommand{\Spp}{\mathrm{Sp}}
\newcommand{\spp}{\mathrm{sp}}
\newcommand{\Ad}{\mathop{\rm Ad}\nolimits}
\newcommand{\ad }{\mathop{\mathrm{ad} }\nolimits}
\newcommand{\dee}{\mathrm{d}}
\newcommand{\comp}{\, \raisebox{2pt}{$\scriptstyle\circ \, $}}
\newcommand{\onehalf}{\mbox{$\frac{\scriptstyle 1}{\scriptstyle 2}\,$}}

\newcommand{\lefthook}{\mbox{$\, \rule{8pt}{.5pt}\rule{.5pt}{6pt}\, \, $}}
\newcommand{\dbydt}{\mbox{${\displaystyle \frac{\dee }{\dee t}}
\rule[-10pt]{.5pt}{25pt} \raisebox{-10pt}{$\, {\scriptstyle t=0}$}$}}

\newcommand{\setrule}{\, \rule[-4pt]{.5pt}{13pt}\, }

\allowdisplaybreaks

\begin{document}
\thispagestyle{empty}
\begin{center}
{\Large \textbf{The momentum mapping of the affine \\ real symplectic group}} \\ 
\mbox{}\vspace{-.2in} \\
\mbox{}\\
Richard Cushman
\mbox{}\\ 
\date{} 
\end{center}
\addtocounter{footnote}{1}
\footnotetext{printed: \today }

In this paper we explain how the cocycle of the momentum map of the action of the 
affine symplectic group on ${\R }^{2n}$ gives rise to a coadjoint orbit of the odd real symplectic 
group with a modulus.   \medskip 

\section{Basic setup}

Consider the set $\mathrm{Af}\Spp ({\R }^{2n}, \omega )$ of invertible affine real 
symplectic mappings 
\begin{displaymath}
(A,a): ({\R }^{2n}, \omega ) \rightarrow  ({\R }^{2n}, \omega ): v \mapsto Av +a . 
\end{displaymath}
Using the multiplication $(A,a) \cdot (B,b) = (AB, Ab+a)$, which corresponds to composition of 
affine linear mappings, $\mathrm{Af}\Spp ({\R }^{2n}, \omega )$ is a group. Identifying $(A,a)$ with the matrix 
\raisebox{2pt}{\tiny $\begin{pmatrix} A & a \\ 0 & 1 \end{pmatrix}$}, 
$\mathrm{Af}\Spp ({\R }^{2n}, \omega )$ becomes 
a closed subgroup of $\Spp ({\R }^{2n}, \omega ) \times {\R }^{2n}$. Thus $\mathrm{Af}\Spp ({\R }^{2n}, \omega )$ is a Lie group. Its Lie algebra is $\mathrm{af}\spp ({\R }^{2n}, \omega ) = 
\{ (X,x) \in \spp ({\R }^{2n}, \omega ) \times {\R }^{2n} \} $ with Lie bracket 
\begin{equation}
[(X,x), (Y,y)] = ([X,Y], Xy-Yx) .
\label{eq-zero}
\end{equation} 

$\mathrm{Af}\Spp ({\R }^{2n}, \omega )$ acts on $({\R }^{2n}, \omega )$ by 
\begin{equation}
\Phi : \mathrm{Af}\Spp ({\R }^{2n}, \omega ) \times ({\R }^{2n}, \omega ) \rightarrow ({\R }^{2n}, \omega ): 
\big( (A,a), v \big) \mapsto Av +a. 
\label{eq-zeronw}
\end{equation}
Since the symplectic form $\omega $ on ${\R }^{2n}$ is invariant under translation, for \linebreak 
every $(A,a) \in \mathrm{Af}\Spp ({\R }^{2n}, \omega )$ the affine linear mapping ${\Phi }_{(A,a)}$ preserves 
$\omega $. The infinitesimal generator $X^{(X,x)}$ of the action $\Phi$ in the 
direction $(X,x) \in \mathrm{af}\spp ({\R }^{2n}, \omega)$ is the vector field 
$X^{(X,x)}(v) = Xv +x$ on ${\R }^{2n}$. We now show \medskip 

\noindent \textbf{Claim 1.1.} The $\mathrm{Af}\Spp ({\R }^{2n}, \omega )$ action $\Phi $ (\ref{eq-zeronw}) is Hamiltonian. \medskip

\noindent \textbf{Proof.} For every $(Y,y) \in \mathrm{af}\spp ({\R }^{2n}, \omega)$ let 
\begin{equation}
J^{(Y,y)}: {\R }^{2n} \rightarrow \R: v \mapsto J^Y(v) + {\omega }^{\sharp}(y)v = 
\onehalf \omega (Yv,v) + \omega (y,v)
\label{eq-zerovn*}
\end{equation}
and set 
\begin{equation}
J: {\R }^{2n} \rightarrow \mathrm{af}\spp ({\R }^{2n}, \omega )^{\ast } ,  
\label{eq-onenw}
\end{equation}
where $J(v)(Y,y) = J^{(Y,y)}(v)$. Then for every $(X,x) \in \mathrm{af}\spp ({\R }^{2n}, \omega)$, 
every $v \in {\R }^{2n}$, and every $w \in T_v{\R }^{2n} = {\R }^{2n}$ 
\begin{align*}
\dee J^{(X,x)}(v)w & = \big( T_vJ\, (X,x) \big) w = \omega (Xv,w) +\omega (x,w) \\
& = \omega (Xv+x, w) = \omega (X^{(X,x)}(v), w),  
\end{align*}
that is, $X^{(X,x)} = X_{J^{(X,x)}}$. Hence the action $\Phi $ is Hamiltonian. 
\hfill $\square $ \medskip 

The above argument shows that the map $J$ (\ref{eq-onenw}) is the momentum map of the Hamiltonian 
action $\Phi $ (\ref{eq-zeronw}). The following discussion is motivated by theorem 11.34
of Souriau \cite[p.143]{souriau}.\medskip 

\noindent \textbf{Lemma 1.2.} The mapping 
\begin{equation}
\sigma : \mathrm{Af}\Spp ({\R }^{2n}, \omega ) \rightarrow \mathrm{af}\spp ({\R }^{2n}, \R )^{\ast }: 
g \mapsto J\big( {\Phi }_g(v) \big) -{\Ad }^T_{g^{-1}}J(v) 
\label{eq-nwzero}
\end{equation}
does not depend on $v \in {\R }^{2n}$. \medskip 

\noindent \textbf{Proof.} For each $\eta \in \mathrm{af}\spp ({\R }^{2n}, \R )$ we have 
\begin{align*}
\dee \big( J \comp {\Phi }_g )^{\eta } (v) & = 
T_v{\Phi }_g \, X^{\eta }(v) \lefthook \omega (v) = X^{{\Ad }_g\eta }(v) \lefthook \omega (v) \\
& = \dee J^{{\Ad }_g\eta }(v) = \dee ({\Ad }^T_{g^{-1}}J)^{\eta }(v),
\end{align*}
that is, $\dee \big( J \comp {\Phi }_g - {\Ad }^T_{g^{-1}}J \big) (v) =0$. Since ${\R }^{2n} $ is 
connected the function $J \comp {\Phi }_g - {\Ad }^T_{g^{-1}}J:{\R }^{2n} \rightarrow \R$ is constant. 
\hfill $\square $ \medskip

\noindent \textbf{Corollary 1.2A} For every $g$, $g' \in \mathrm{Af}\Spp ({\R }^{2n}, \omega )$ 
\begin{equation}
\sigma (gg') = \sigma (g) + {\Ad }^T_{g^{-1}}\sigma (g'). 
\label{eq-nwthreestar}
\end{equation}

\noindent \textbf{Proof.} We compute. 
\begin{align}
\sigma (gg') & = J \comp {\Phi }_{gg'} -{\Ad }^T_{(gg')^{-1}}J  \notag \\
& = (J\comp {\Phi }_g -{\Ad }^T_{g^{-1}}J) \comp {\Phi }_{g'} + 
{\Ad }^T_{g^{-1}}(J \comp {\Phi }_{g'} - {\Ad }^T_{(g')^{-1}}J) \notag \\
& = \sigma(g) + {\Ad }^T_{g^{-1}}\sigma (g'). \tag*{$\square $}
\end{align}

Evaluating $\sigma $ (\ref{eq-nwzero}) at $\exp t \eta $ and then at $\zeta \in \mathrm{af}\spp ({\R }^{2n}, \R )$ gives 
\begin{equation}
\big( \sigma (\exp t\eta ) \big) \zeta = J^{\zeta }\big( {\Phi }_{\exp t \eta }(v) \big) - 
\big( {\Ad }^T_{\exp -t\eta } J(v) \big) \zeta . 
\label{eq-nwone}
\end{equation}
Differentiating (\ref{eq-nwone}) with respect to $t$ and then setting $t$ equal to $0$ gives
\begin{align}
(T_e\sigma \, \eta )\zeta & = \dee J^{\zeta }(v) X^{\eta }(v) + \big( {\ad }^T_{\eta }J(v) \big) \zeta \notag \\
& = L_{X^{\eta }}J^{\zeta}(v) + J(v){\ad }_{\eta }\zeta 
= \{ J^{\zeta }, J^{\eta } \} (v) - J^{[\zeta , \eta ]}(v). 
\label{eq-nwtwo}
\end{align}

Let ${\Sigma }^{\sharp}: \mathrm{af}\spp ({\R }^{2n}, \R ) \rightarrow \mathrm{af}\spp ({\R }^{2n}, \R )^{\ast }$ 
be the linear mapping $\eta \mapsto {\Sigma }^{\sharp}(\eta ) = -(T_e\sigma ) \eta \in 
\mathrm{af}\spp ({\R }^{2n}, \R )^{\ast }$. Equation (\ref{eq-nwtwo}) may be written as 
\begin{equation}
\{ J^{\eta }, J^{\zeta } \} (v) = J^{[\eta , \zeta ]}(v) + \Sigma (\eta , \zeta ) , 
\label{eq-nwthree}
\end{equation}
where $\Sigma (\eta , \zeta ) = {\Sigma }^{\sharp}(\eta ) \zeta $. From equation (\ref{eq-nwthree}) 
it follows that the bilinear map $\Sigma $ is skew symmetric.  \medskip

\noindent \textbf{Lemma 1.3.} $\Sigma $ is an $\mathrm{af}\spp ({\R }^{2n}, \R )$ cocycle, that is, 
for every $\xi $, $\eta $, and $\zeta \in \mathrm{af}\spp ({\R }^{2n}, \R )$ 
\begin{equation}
\Sigma (\zeta , [\xi , \eta ]) = \Sigma ([\zeta , \xi ], \eta ) + \Sigma (\xi , [\zeta , \eta ]). 
\label{eq-nwfour}
\end{equation}

\noindent \textbf{Proof.} Since $\big( C^{\infty}({\R }^{2n}), \{ \, \, , \, \, \}  \big) $ is a Lie algebra
\begin{displaymath}
\{ J^{\zeta }, \{ J^{\xi }, J^{\eta } \} \} = \{ \{ J^{\zeta }, J^{\xi } \}, J^{\eta } \} + \{ J^{\xi }, \{ J^{\zeta }, 
J^{\eta } \} \} .
\end{displaymath}
Using equation (\ref{eq-nwthree}) the above equation reads
\begin{align*}
\{ J^{\zeta }, J^{[\xi , \eta ]} \} + \{ J^{\zeta }, \Sigma (\xi ,\eta ) \} & = 
\{ J^{[\zeta , \xi ]}, J^{\eta } \} + \{ \Sigma (\zeta , \xi ), J^{\eta } \} \\
& \hspace{-1in} + \{ J^{\xi }, J^{[\zeta , \eta ]} \} + \{ J^{\zeta }, \Sigma (\zeta ,\eta ) \} .
\end{align*}
Using (\ref{eq-nwthree}) again gives 
\begin{displaymath}
J^{[\zeta , [\xi , \eta ]} + \Sigma (\zeta , [\xi , \eta ]) = J^{[[\zeta, \xi ], \eta ]} + \Sigma ([\zeta , \xi ], \eta ) 
+J^{[\xi , [\zeta , \eta ]} + \Sigma (\xi , [\zeta , \eta ]), 
\end{displaymath}
since $\Sigma (\xi , \eta )$, $\Sigma (\zeta ,\xi )$ and $\Sigma (\zeta , \eta )$ are constant functions on 
${\R }^{2n}$. By linearity and the Jacobi identity $J^{[\zeta ,[\xi \eta ]]} = J^{[[\zeta ,\xi ] ,\eta ]} 
+J^{[\xi , [\zeta , \eta ]]}$ equation (\ref{eq-nwfour}) holds. \hfill $\square $ \medskip 

\noindent \textbf{Claim 1.4.} The momentum map $J$ (\ref{eq-onenw}) has the cocycle 
\begin{equation}
\Sigma : \mathrm{af}\spp ({\R }^{2n}, \omega) \times \mathrm{af}\spp ({\R }^{2n}, \omega )  \rightarrow 
\R :   \big( (Y,y) , (Z,z) \big) \mapsto {\omega }(y,z) . 
\label{eq-twonw}
\end{equation}

\noindent \textbf{Proof.} We compute. 
\begin{align*}
\{ J^{(Y,y)}, J^{(Z,z)} \} (v) & = \big( L_{X^{(Z,z)}} J^{(Y,y)}\big) (v) = 
\dee J^{(Y,y)} (v) X^{(Z,z)}(v)   \\
& = \omega (Yv, Zv+z) + \omega (y, Zv+z)  \\
& = \omega (Yv, Zv) + \omega (Yv,z) + \omega (y,Zv) + \omega (y,z) .
\end{align*}
Now
\begin{align*} 
\onehalf \omega ([Y,Z]v,v)  & = \onehalf \omega \big( (YZ -ZY)v,v) 
= \onehalf \omega (YZv,v) -\onehalf \omega (ZYv,v) \\
& = -\onehalf \omega (Zv,Yv) + \onehalf \omega (Yv,Zv)  = \omega (Yv,Zv). 
\end{align*}
Thus 
\begin{align}
\{ J^{(Y,y)}, J^{(Z,z)} \} (v) & = \onehalf \omega ([Y,Z]v,v) + \omega (Yz -Zy, v) + \omega (y,z) \notag \\
& = J^{[Y,Z]}(v) + \omega (Yz-Zy, v) + \omega (y,z) \notag \\ 
& = J^{[(Y,y), (Z,z)]}(v)  + \omega (y,z) . \tag*{$\square $}
\end{align}

Define the map 
\begin{equation}
\begin{array}{l}
\Psi : \mathrm{Af}\Spp ({\R}^{2n}, \omega ) \times \mathrm{af}\spp ({\R }^{2n}, \omega )^{\ast } 
\rightarrow \mathrm{af}\spp ({\R }^{2n}, \omega )^{\ast } \\
\hspace{1in}(g, \alpha ) \longmapsto {\Ad }^T_{g^{-1}}\alpha + \sigma (g) , 
\label{eq-onevnw**}
\end{array}
\end{equation} 
where $\sigma $ is given by equation (\ref{eq-nwzero}). \medskip 

\noindent \textbf{Claim 1.5} The map $\Psi $ (\ref{eq-onevnw**}) is an action of 
$\mathrm{Af}\Spp ({\R}^{2n}, \omega )$ on $\mathrm{af}\spp ({\R }^{2n}, \omega )^{\ast }$. \medskip 

\noindent \textbf{Proof.} We compute. For $g$, $g' \in \mathrm{Af}\Spp ({\R}^{2n}, \omega )$ and $\alpha \in \mathrm{af}\spp ({\R }^{2n}, \omega )^{\ast }$ we have 
\begin{align}
{\Psi }_{gg'} \alpha & = {\Ad }^T_{(gg')^{-1}}\alpha + \sigma (gg') \notag \\
& = {\Ad }^T_{g^{-1}}({\Ad }^T_{(g')^{-1}}\alpha ) + \sigma (g) + {\Ad }^T_{g^{-1}}\sigma (g') \notag \\
& \hspace{.75in}\mbox{using corollary 1.2A} \notag \\
& = {\Ad }^T_{g^{-1}}\big( {\Ad }^T_{(g')^{-1}}\alpha + \sigma (g') \big) + \sigma (g) \notag \\
& = {\Ad }^T_{g^{-1}}({\Psi }_{g'}\alpha ) + \sigma (g) = {\Psi }_g({\Psi }_{g'}\alpha ). 
\tag*{$\square $} 
\end{align}

\noindent \textbf{Claim 1.6} The momentum mapping $J$ (\ref{eq-onenw}) is coadjoint 
equivariant under the action $\Psi $ (\ref{eq-onevnw**}). \medskip 

\noindent \textbf{Proof.} We compute. For every $g \in \mathrm{Af}\Spp ({\R}^{2n}, \omega )$ and 
every $w \in {\R }^{2n}$ 
\begin{align*}
{\Psi }_g\big( J(w) \big) & = {\Ad }^T_{g^{-1}}J(w) + \sigma (g), \, \, \, \mbox{using (\ref{eq-onevnw**})} \\
& = {\Ad }^T_{g^{-1}}J(w) + J\big( {\Phi }_g(w) \big) - {\Ad }^T_{g^{-1}}J(w), \, \, 
\mbox{using (\ref{eq-nwzero}) } \\
& = J\big( {\Phi }_g(w) \big) . \tag*{$\square $}
\end{align*}

\section{Extension}

Following Wallach \cite{wallach} we find a central extension of Lie algebra $\mathrm{af}\spp ({\R }^{2n}, \omega )$ the dual of whose adjoint map is 
\clearpage
\begin{align} 
T_e\Psi (X,x) \alpha & = \dbydt {\Psi }_{\exp t(X,x)}\alpha 
= -{\ad }^T_{(X,x)}\alpha + T_e\sigma (X,x) \notag \\
& = -{\ad }^T_{(X,x)}\alpha + {\Sigma }^{\sharp}(X,x), 
\label{eq-twovnw}
\end{align}
where $(X,x) \in \spp ({\R }^{2n}, \omega ) \times {\R }^{2n} = \mathrm{af}\spp ({\R }^{2n}, \omega )$. \medskip 

Consider the Lie algebra $\widehat{\mathfrak{g}} = \{ (X,v, \xi ) \in \mathrm{af}\spp ({\R }^{2n}, \omega ) \times \R \}$ whose Lie bracket is 
\begin{equation}
[(X,v, \xi ), (Y,w, \eta )] = ([X,Y], Xw-Yv, \omega (v,w)). 
\label{eq-onevnw*}
\end{equation}
From (\ref{eq-onevnw*}) it follows that $(0,0,1)$ lies in the center of $\widehat{\mathfrak{g}}$, that is, 
$[(0,0,1), (X,v,\xi )]$ $= (0,0,0)$. Also 
\begin{equation}
[(X,v,0), (Y,v,0)] = ([X,Y], Xw-Yv,0) + \omega (v,w)\, (0,0,1). 
\label{eq-fourvnw}
\end{equation}
Thus the Lie algebra $\widehat{\mathfrak{g}}$ is a central extension of the Lie algebra 
$\mathrm{af}\spp ({\R }^{2n}, \omega )$ by the cocycle $\omega $. Since we can write 
(\ref{eq-onevnw*}) as 
\begin{equation}
{\ad }_{[X,x, \xi ]}[Y,y, \eta ] = {\ad }_{[X,x]}[Y,y] + \Sigma (\xi ,\eta ), 
\label{eq-twovnw*}
\end{equation}
$\widehat{\mathfrak{g}}$ is the sought for Lie algebra. 

\section{The odd real symplectic group}

We now find a connected linear Lie group $\widehat{G}$ whose Lie algebra is $\widehat{\mathfrak{g}}$. 
Consider the group $\widehat{G} \subseteq \mathrm{Af}\Spp ({\R }^{2n}, \onehalf \omega ) \times \R$ with 
multiplication 
\begin{displaymath}
\big( (A,v),r \big) \cdot \big( (B,w),s \big) = \big( (AB, Aw+v), r+s+ \onehalf \omega (A^{-1}v,w ) \big) . 
\end{displaymath}
The map 
\begin{displaymath}
\widehat{\pi }: \widehat{G} \rightarrow \mathrm{Af}\Spp ({\R }^{2n}, \omega ): (A,v,r) \mapsto (A,v)
\end{displaymath}
is a surjective group homomorphism, whose kernel is the normal subgroup 
$\widehat{Z} = \{ (I,0,r) \in \widehat{G} \setrule \, r \in \R \}$, which is the center of 
$\widehat{G}$. Note that ${\pi }_1(\widehat{G}) = 
{\pi }_1(\Spp ({\R }^{2n}, \omega )) = \Z$. $\widehat{G}$ is a Lie group with Lie algebra 
$\widehat{\mathfrak{g}}$, whose Lie bracket is given by (\ref{eq-onevnw*}). \medskip 

\noindent \textbf{Claim 3.1.} The group $\widehat{G}$ is isomorphic to the odd real symplectic group. \medskip

\noindent \textbf{Proof.} Consider the map 
\begin{equation}
\rho : \widehat{G} \rightarrow \Gl ({\R }^{2n+2}, \R ): (A,v,r) \mapsto 
\mbox{\footnotesize $\begin{pmatrix} 
1 & 0 & 0 \\
v & A & 0 \\
r & \onehalf {\omega }^{\sharp}(A^{-1}v) & 1 
\end{pmatrix} .$}
\label{eq-fivevnw*}
\end{equation}
The map $\rho $ is an injective homomorphism. Here we have 
$\onehalf {\omega }^{\sharp}(A^{-1}v) =$ \newline  
$-\onehalf (v^TJA)$, since 
for every $z \in {\R }^{2n}$
\begin{align*}
\onehalf {\omega }^{\sharp}(A^{-1}v)z & = \onehalf \omega (A^{-1}v,z) = \onehalf \omega (v, Az), \, \, 
\mbox{because $A \in \Spp ({\R }^{2n}, \onehalf \omega )$} \\
& = -\onehalf \omega (Az, v) =  -\onehalf (v^T J A)z.
\end{align*} 
The following calculation shows that $\rho $ is a homomorphism.  
\begin{align}
\rho \big( (A,v,r) \cdot (B,w,s) \big) & = \rho \big( AB, Aw +v, r+s + \onehalf \omega (A^{-1}v,w) \big) \notag \\
& \hspace{-.5in} = \mbox{\footnotesize $\begin{pmatrix}
1 & 0 & 0 \\
v+Aw & AB & 0 \\
r+s + \onehalf {\omega }^{\sharp}(A^{-1}v)w & \onehalf {\omega}^{\sharp}\big( (AB)^{-1}(v + Aw) \big) & 1 
\end{pmatrix}$.} 
\label{eq-sixvnw**}
\end{align}
Now 
\begin{align}
\onehalf {\omega}^{\sharp}\big( (AB)^{-1}(v + Aw) \big) & = \onehalf {\omega }^{\sharp}\big( B^{-1}(A^{-1}v) \big) + \onehalf {\omega }^{\sharp}(B^{-1}w) \notag \\
&\hspace{-1in} = B^T \onehalf {\omega }^{\sharp}(A^{-1}v) + \onehalf {\omega }^{\sharp}(B^{-1}w), \, \, \, 
\mbox{since $B \in \Spp ({\R }^{2n}, \onehalf \omega )$} \notag \\
&\hspace{-1in} = \onehalf {\omega }^{\sharp}(A^{-1}v)B + \onehalf {\omega }^{\sharp}(B^{-1}w) . 
\label{eq-fivevnw}
\end{align} 
But 
\begin{align}
\rho (A,v,r) \, \rho (B,w,s) & = 
\mbox{\footnotesize $\begin{pmatrix} 
1 & 0 & 0 \\
v & A & 0 \\
r &\onehalf {\omega }^{\sharp}(A^{-1}v) & 1 \end{pmatrix} \, 
\begin{pmatrix} 
1 & 0 & 0 \\
w & B & 0 \\
s &\onehalf {\omega }^{\sharp}(B^{-1}w) & 1 \end{pmatrix}$} \notag \\
&\hspace{-.5in} = \mbox{\footnotesize $\begin{pmatrix}
1 & 0 & 0 \\
v+Aw & AB & 0 \\
r+s + \onehalf {\omega }^{\sharp}(A^{-1}v)w & \onehalf {\omega }^{\sharp}(A^{-1}v)B + 
\onehalf {\omega }^{\sharp}(B^{-1}w) & 1 
\end{pmatrix}$.} 
\label{eq-sixvnw*}
\end{align}
Using (\ref{eq-fivevnw}) we see that the right hand sides of equations (\ref{eq-sixvnw**}) and 
(\ref{eq-sixvnw*}) are equal, that is, 
\begin{displaymath}
\rho \big( (A,v,r) \cdot (B,w,s) \big) = \rho (A,v,r) \, \rho (B,w,s). 
\end{displaymath}
Thus the map $\rho $ (\ref{eq-fivevnw*}) is a group homomorphism. The map $\rho $ is injective, for if $\rho (A,v,r) = (I_{2n}, 0,0)$, then $v=0$ and $r=0$. So 
$(A,v,r) = (I_{2n},0,0)$, the  identity element of $\widehat{G}$. \medskip 

Since $\widehat{G}$ is a subgroup of $\mathrm{Af}\Spp ({\R }^{2n}, \onehalf \omega ) \times \R$, it follows that 
if $(A,v,r) \in \widehat{G}$, then $A \in \Spp ({\R }^{2n}, \onehalf \omega  )$. Thus the image of the map 
$\rho $ is contained in $\Spp ({\R }^{2n+2}, \Omega )$. Here the matrix of the symplectic form 
$\Omega $ with respect to the basis $\{ e_0, e_1, \ldots , e_n, f_1, \ldots , f_n, f_{n+1} \}$ of 
${\R }^{2n+2}$ is {\tiny $\begin{pmatrix} 0 & 0 &1 \\
0 & \onehalf J & 0 \\ -1 & 0 & 0 \end{pmatrix} $,} and $J=$ \linebreak 
{\tiny $\begin{pmatrix} 
0 & I_n \\ -I_n & 0 \end{pmatrix}$} is the matrix of the symplectic form $\omega $ with respect to 
the basis $\{ e_1, \ldots , e_n, f_1, $ $ \ldots , f_n \} $ of ${\R }^{2n}$.  The 
image of the map $\rho $ (\ref{eq-fivevnw*}) is the \emph{odd real symplectic group} 
${\Spp ({\R}^{2n+2}, \Omega )}_{f_{n+1}} = 
\{ \mathcal{A} \in \Spp ({\R }^{2n+2}, \Omega ) \setrule \, \mathcal{A}f_{n+1} = f_{n+1} \}$. Consequently, 
$\widehat{G}$ is isomorphic to ${\Spp ({\R}^{2n+2}, \Omega )}_{f_{n+1}}$. \hfill $\square $ \medskip  

The Lie algebra ${\spp ({\R }^{2n+2}, \Omega )}_{f_{n+1}}$ of the Lie group 
${\Spp ({\R}^{2n+2}, \Omega )}_{f_{n+1}}$ is 
\begin{displaymath}
\{ \widehat{X} = \mbox{\tiny $\begin{pmatrix} 0 & 0 & 0 \\ x & X & 0 \\ \xi & \onehalf {\omega }^{\sharp}(x) & 0 
\end{pmatrix}$} \in \Spp ({\R }^{2n+2}, \Omega ) \setrule  x \in {\R }^{2n}, \, 
X \in \spp ({\R }^{2n}, \onehalf \omega ), \, \mathrm{and}\, \, \xi \in \R \} 
\end{displaymath}
with Lie bracket 
\begin{align}
[\widehat{X}, \widehat{Y}] & = \left[ \mbox{\tiny $\begin{pmatrix} 
0 & 0 & 0 \\ x & X & 0 \\ \xi & \onehalf {\omega }^{\sharp}(x) & 0 \end{pmatrix}, \, 
\begin{pmatrix} 0 & 0 & 0 \\ y & Y & 0 \\ \eta  & \onehalf {\omega }^{\sharp}(y) & 0 \end{pmatrix}$} \right] 
= \mbox{\tiny $\begin{pmatrix} 0 & 0 & 0 \\ Xy-Yx & [X,Y] & 0 \\ \omega (x,y) & 
\onehalf {\omega }^{\sharp}(Xy-Yx) & 0 \end{pmatrix} $.} 
\label{eq-eight}
\end{align}
Here $\onehalf {\omega }^{\sharp}(x) = -\onehalf x^TJ$, since for each $z \in {\R }^{2n}$ 
\begin{align*}
\onehalf {\omega }^{\sharp}(x)z & = \onehalf \omega (x,z) = - \onehalf \omega (z,x) = 
(-\onehalf x^TJ)z . 
\end{align*}
The map 
\begin{displaymath}
\mu :\widehat{\mathfrak{g}} \rightarrow {\spp ({\R }^{2n+2}, \Omega )}_{f_{n+1}}: 
(X,x, \xi ) \mapsto \widehat{X} = \mbox{\footnotesize $\begin{pmatrix} 
0 & 0 & 0 \\ x & X & 0 \\ \xi & \onehalf {\omega }^{\sharp}(x) & 0 \end{pmatrix} $} 
\end{displaymath}
is a Lie algebra isomorphism, because it is a bijective linear map and 
\begin{align*}
\mu \big( [(X,x,\xi ), (Y,y,\eta ) ] \big) & = \mu \big( [X,Y], Xy-Yx, \omega (x,y)  \big) \\
&\hspace{-.75in} = \mbox{\tiny $\begin{pmatrix} 0 & 0 & 0 \\ Xy-Yx & [X,Y] & 0 \\ \omega (x,y) & 
\onehalf {\omega }^{\sharp}(Xy-Yx) & 0 \end{pmatrix} $} 
=[\widehat{X}, \widehat{Y}] = [\mu (X,x, \xi), \mu (Y,y, \eta ) ] . 
\end{align*}
This verifies that the Lie algebra of the Lie group $\widehat{G}$ has Lie bracket given by 
(\ref{eq-onevnw*}), because the group $\widehat{G}$ is isomorphic to 
${\Spp ({\R }^{2n+2}, \Omega )}_{f_{n+1}}$. \medskip 

\noindent \textbf{Claim 3.2.} The action 
\begin{equation}
\begin{array}{l}
\widehat{\Phi }: {\Spp ({\R }^{2n+2}, \Omega )}_{f_{n+1}} \times ({\R }^{2n}, \omega ) \rightarrow 
({\R }^{2n}, \omega ) :  \\ 
\rule{0pt}{22pt}\hspace{.25in} \Big( \mbox{\footnotesize $\begin{pmatrix} 
1 & 0 & 0 \\ v & A & 0 \\ r & \onehalf {\omega }^{\sharp}(A^{-1}v) & 1 \end{pmatrix} $}, \, 
w  \Big) \mapsto  Aw +v 
\end{array} 
\label{eq-nine}
\end{equation}
is Hamiltonian. \medskip 

\noindent \textbf{Proof.} The infinitesimal generator $X^{\widehat{X}}$ of the action 
$\widehat{\Phi }$ in the direction $\widehat{X} \in {\spp ({\R }^{2n+2}, \Omega )}_{f_{n+1}}$ is the vector field 
$X^{\widehat{X}}(w) = X(w) +x$ on ${\R }^{2n}$. For every $\widehat{Y}=${\tiny $\begin{pmatrix} 
0 & 0 & 0 \\ y & Y & 0 \\ \eta & \onehalf {\omega }^{\sharp}(y) & 0 \end{pmatrix} $}
$\in {\spp ({\R }^{2n+2},\Omega )}_{f_{n+1}}$ let 
\begin{equation}
J^{\, \widehat{Y}}: {\R }^{2n} \rightarrow \R : w \mapsto 
\onehalf \omega (Yw, w) + \omega (y,w) + \eta . 
\label{eq-ten}
\end{equation}
Then 
\begin{displaymath}
\dee {\widehat{J}}^{\, \, \widehat{Y}}(v)w  = T_v\widehat{J}(\, \widehat{Y})w = 
\omega (Yv,w) + \omega (y,w) = \omega \big( X^{\widehat{Y}}(v), w). 
\end{displaymath}
Thus $X^{\widehat{Y}} = X_{J^{\widehat{Y}}}$. So the action $\widehat{\Phi }$ (\ref{eq-nine}) is Hamiltonian. 
Since the mapping $\widehat{Y} \mapsto J^{\widehat{Y}}(w)$  is linear for every $w \in {\R }^{2n}$, the action 
$\widehat{\Phi }$ (\ref{eq-nine}) has a momentum mapping  
\begin{equation}
\widehat{J}: {\R }^{2n} \rightarrow {\spp ({\R }^{2n+2}, \Omega )}^{\ast}_{f_{n+1}} , 
\label{eq-eleven*}
\end{equation}
with $\widehat{J}(w) \widehat{Y} = {\widehat{J}}^{\, \, \widehat{Y}}(w)$. \hfill $\square $ \medskip 

\noindent \textbf{Claim 3.3.} The momentum mapping $\widehat{J}$ (\ref{eq-eleven*}) of the 
${\Spp ({\R }^{2n+2}, \Omega )}_{f_{n+1}}$ action $\widehat{\Phi }$ (\ref{eq-nine}) is 
coadjoint equivariant, that is, 
\begin{equation}
\widehat{J}\big( {\widehat{\Phi }}_g(w) \big) = {\Ad }^T_{g^{-1}}\widehat{J}(w) 
\label{eq-eleven}
\end{equation}
for every $g \in {\Spp ({\R }^{2n+2}, \Omega )}_{f_{n+1}}$ and every $w \in {\R }^{2n}$. \medskip 

\noindent \textbf{Proof.} It is enough to show that 
\begin{equation}
\{ {\widehat{J}}^{\, \, \widehat{Y}}, {\widehat{J}}^{\, \, \widehat{Z}} \} = 
{\widehat{J}}^{\, \, [\widehat{Y}, \widehat{Z}]}, 
\, \, \mbox{for every $\widehat{Y}$, $\widehat{Z} \in \widehat{\mathfrak{g}}$,} 
\label{eq-twelve*}
\end{equation}
because (\ref{eq-twelve*}) is the infinitesimalization of the coadjoint equivariance condition (\ref{eq-eleven}) and ${\Spp ({\R }^{2n+2}, \Omega )}_{f_{n+1}}$ is generated by elements which lie in the image 
of the exponential mapping, since it is connected. We compute 
\begin{align}
{\widehat{J}}^{\, \, [\widehat{Y}, \widehat{Z}]}(w) & = 
\onehalf \omega ([Y,Z]w,w) + \omega (Yz-Zy, w) +\omega (y,z), \notag \\ 
& \hspace{.5in} \mbox{using equations (\ref{eq-eight}) and (\ref{eq-ten})} \notag 
\\ 
& = \omega (Yw, Zw) + \omega (Yw,z) +\omega (y, Zw) + \omega (y,z) \notag 
\\
& = \omega (Yw, Zw+z) + \omega (y, Zw+z) 
= \dee {\widehat{J}}^{\, \, \widehat{Y}}(w) X^{\widehat{Z}}(w) \notag 
\\
& =\big(  L_{X_{{\widehat{J}}^{\, \widehat{Z}}}}{\widehat{J}}^{\, \, \widehat{Y}}\big) (w) = 
 \{ {\widehat{J}}^{\, \, \widehat{Y}} \! , {\widehat{J}}^{\, \widehat{Z}} \} (w). \tag*{$\square $}
\end{align}

\section{Coadjoint orbit}

In this section using results of \cite{cushman22} we algebraically classify the coadjoint 
orbit $\mathcal{O}\big(\widehat{J}(e_1) \big)$ of the odd real symplectic group 
${\Spp ({\R }^{2n+2}, \Omega )}_{f_{n+1}}$ through $\widehat{J}(e_1) \in 
{\spp ({\R }^{2n+2}, \Omega )}^{\ast}_{f_{n+1}}$. We show that this coadjoint orbit has a modulus. \medskip 

\vspace{-.15in} First we note that the action $\widehat{\Phi }$ (\ref{eq-nine}) of
${\Spp ({\R }^{2n+2}, \Omega )}_{f_{n+1}}$ on ${\R }^{2n}$ is 
transitive. Thus to determine the ${\Spp ({\R }^{2n+2}, \Omega )}_{f_{n+1}}$ coadjoint orbit 
through $\widehat{J}(w)$ for a fixed $w \in {\R }^{2n}$, it suffices to determine the coadjoint orbit 
$\mathcal{O}\big(\widehat{J}(e_1) \big)$ through $\widehat{J}(e_1) \in 
{\spp ({\R }^{2n+2}, \Omega )}^{\ast }_{f_{n+1}}$. We have 
\begin{align*}
\widehat{Y} & = \mbox{\footnotesize $\left( \begin{array}{c|c|c|c}
0 & 0 & 0 & 0 \\ \hline
\rule{0pt}{9pt}y^1 & A & B & 0 \\ \hline
\rule{0pt}{10pt} y^2 & C & -A^T & 0 \\ \hline  
\rule{0pt}{10pt} \eta & \onehalf (y^2)^T & -\onehalf (y^1)^T & 0 \end{array} \right) $ ,}
\end{align*}
where $Y = (Y_{ij}) =${\tiny $\begin{pmatrix} A & B \\ C & -A^T \end{pmatrix}$}$\in \spp ({\R }^{2n}, 
\onehalf \omega )$ 
with $A$, $B$, and $C \in \gl ({\R }^n, \R)$, where $B=B^T$ and 
$C = C^T$; $y^T = \big( (y^1)^T \! \mid \! (y^2)^T \big) = 
\big( y_1, \ldots , y_n \! \mid \! y_{n+1} , \ldots , y_{2n} \big) $ 
$\in {\R }^{2n}$; and $\eta \in \R $. Then using (\ref{eq-ten}) we get
\begin{align*}
\widehat{J}(e_1)\widehat{Y} & = \onehalf e^T_1\mbox{\footnotesize $\begin{pmatrix} C & -A^T \\ -A & - B \end{pmatrix}$}e_1 + e^T_1 \mbox{\footnotesize $\begin{pmatrix} y^2 \\ -y^1 \end{pmatrix}$}  + \eta \\
& = \onehalf Y_{n+1, 1} + y_{n+1} +  \eta . 
\end{align*}
Therefore 
\begin{equation}
\widehat{J}(e_1) = \onehalf E^{\ast }_{n+1, 1} + \onehalf E^{\ast}_{n+1,0}  + 
E^{\ast }_{2n+1,1} + E^{\ast }_{2n+1, 0}. 
\label{eq-thirteen}
\end{equation}
With respect to the basis ${\{ E^{\ast }_{ij} \} }^{2n+1}_{i,j = 0}$, which is dual to 
the standard basis ${\{ E_{ij} \} }^{2n+1}_{i,j=0}$ of $\gl ({\R }^{2n+2}, \R )$ where $E_{ij} = 
({\delta }_{ik}{\delta }_{j\ell})$, we have 
\begin{equation}
\widehat{J}(e_1) = 
\mbox{\footnotesize $\left( \begin{array}{c|c|c|c} 
0 & 0 & 0 & 0 \\ \hline 
0 & 0 & 0 & 0 \\ \hline 
r & D & 0 & 0 \\ \hline 
\rule{0pt}{10pt} 1 & 2 r^T& 0 & 0 \end{array} \right) $} 
\in {\spp ({\R }^{2n+2}, \Omega )}^{\ast }_{f_{n+1}}. 
\label{eq-fourteen}
\end{equation}
Here $r^T = (1/2, 0, \ldots , 0) \in {\R }^n$ and $D= \mathrm{diag}(1/2, 0, \ldots , 0) \in \gl ({\R }^n, \R )$. \medskip 

We now use results of \cite{cushman22} to algebraically characterize the coadjoint orbit 
$\mathcal{O}\big( \widehat{J}(e_1) \big)$. The \emph{affine cotype} $\nabla$ represented by the tuple 
$({\R }^{2n+2}, Z, f_{n+1}; \Omega )$ corresponds to the coadjoint orbit 
$\mathcal{O}\big( \widehat{J}(e_1) \big)$, see proposition 4 of \cite{cushman22}. Here 
\begin{displaymath}
Z = \widehat{J}(e_1)^T = \mbox{\footnotesize $\left( \begin{array}{c|c|c|c} 
0 & 0 & r^T & 1 \\ \hline
0 & 0 & D & 2r \\ \hline 
0 & 0 & 0 & 0 \\ \hline 
0 & 0 & 0 & 0 \end{array} \right) $} \in \spp ({\R }^{2n+2}, \Omega ) . 
\end{displaymath}
Since $Z^2 =0$, the cotype $\nabla $ is nilpotent and has height $1$. The parameter of 
$\nabla $ is $1$. Thus by proposition 7 of \cite{cushman22} the affine cotype 
$\nabla = {\nabla }_1(0), 1 + \Delta $. 
Here ${\nabla }_1(0), 1$ is the nilpotent indecomposable cotype of height $1$ and 
\emph{modulus} $1$, which is represented by the tuple 
$\big( {\R }^2, $ {\tiny $\begin{pmatrix} 0 & 1 \\ 0 & 0 \end{pmatrix} $}, $f_1;$ 
{\tiny $\begin{pmatrix} 0 & 1 \\ -1 & 0 \end{pmatrix} $} $\big) $, and $\Delta $ is a semisimple type, 
see \cite{burgoyne-cushman}. Since $\Delta $ is a nilpotent of height $0$, it equals $0$.  \medskip 

We give another argument which proves the above assertion. The tuples $({\R }^{2n+2}, Z, f_{n+1}; \Omega )$ and $({\R }^{2n+2}, PZP^{-1},f_{n+1}; \Omega )$ are equivalent, when $P \in \Spp ({\R }^{2n+2}, \Omega )_{f_{n+1}}$ and thus represent the same cotype $\nabla $. Let $P=${\tiny $\begin{pmatrix} 
1 & 0 & 0 \\ d & I_{2n} & 0 \\ \rule{0pt}{7pt}f & -\onehalf d^T J  & 1 \end{pmatrix}$} 
$\in \Spp ({\R }^{2n+2}, \Omega )_{f_{n+1}}$. We compute.
\begin{align*}
PZP^{-1} & = \mbox{\footnotesize $\begin{pmatrix} 1 & 0 & 0 \\ d & I_{2n} & 0 \\
f & -\onehalf d^TJ & 1\end{pmatrix} \, 
\begin{pmatrix} 0 & {\widetilde{r}}^{\, T} & 1 \\ 0 & \widetilde{D} & \widetilde{s} \\ 0 & 0 & 0 \end{pmatrix} \, 
\begin{pmatrix} 1 & 0 & 0 \\ -d & I_{2n} & 0 \\ -f & \onehalf d^TJ & 1 \end{pmatrix} $, } \\
& \hspace{.5in}\parbox{3.5in}{where ${\widetilde{r}}^{\, T} = ( 0 \! \mid \! r^T )$, 
${\widetilde{s}}^{\, T}= (2J\, \widetilde{r})^T = (2r \! \mid \! 0 )$, and $\widetilde{D}=
${\tiny $\left( \begin{array}{c|c} 0 & D \\ \hline  0 & 0 \end{array} \right) $} } \\
& = \mbox{\footnotesize $\begin{pmatrix} 1 & 0 & 0 \\ d & I_{2n} & 0 \\
f & -\onehalf d^TJ & 1\end{pmatrix} \, \begin{pmatrix} -{\widetilde{r}}^{\, T}d -f & 
{\widetilde{r}}^{\, T} + \onehalf d^TJ & 1 \\ -\widetilde{D}d - f \widetilde{s} & 
\widetilde{D} + \onehalf \widetilde{s} \otimes d^TJ & \widetilde{s} \\
0 & 0 & 0 \end{pmatrix} $.} 
\end{align*}
Choose $d = - \widetilde{s}$ and set $f=0$. Then 
\begin{align*}
&-r^Td -f = {\widetilde{r}}^{\, T}\widetilde{s} = (0 \! \mid \! r^T) \, \mbox{\footnotesize $\left( \begin{array}{c} 
2r \\ \hline 0 \end{array} \right) $} = 0 \\ 
&{\widetilde{r}}^{\, T} + \onehalf d^T J = (0\! \mid \! r^T) + \onehalf (-2r^T \! \mid \! 0) \mbox{\footnotesize 
$\begin{pmatrix} 0 & I \\ -I & 0 \end{pmatrix}$} = (0\! \mid \! r^T) + (0\! \mid \! -r^T) =0 \\
&-\widetilde{D}d - f \, \widetilde{s} = -\widetilde{D}\widetilde{s} = \mbox{\footnotesize $\left( \begin{array}{c|c} 
0 & D \\ \hline 0 & 0 \end{array} \right) \, \left( \begin{array}{c} 2r \\ \hline 0 \end{array} \right) $} =0 \\
&\widetilde{D} + \onehalf \widetilde{s} \otimes d^TJ = \widetilde{D} + \onehalf \widetilde{s} \otimes 
( 0 \! \mid \! -2r) = \widetilde{D} - \onehalf e_1 \otimes f^T_1 =0. 
\end{align*}
Therefore 
\begin{align*}
PZP^{-1} & = \mbox{\footnotesize $\begin{pmatrix} 1 & 0 & 0 \\ -\widetilde{s} & I_{2n} & 0 \\ 0 & 
\onehalf {\widetilde{s}}^{\, T}J & 1 \end{pmatrix} \, \begin{pmatrix} 0 & 0 & 1 \\ 0 & 0 & \widetilde{s} \\ 
0 & 0 & 0 \end{pmatrix} $} = \mbox{\footnotesize $\begin{pmatrix} 0 & 0 & 1 \\ 0 & 0 & 0 \\ 0 & 0 & 0 
\end{pmatrix} $.}
\end{align*}

\end{document}